\documentclass[final,1p,times]{elsarticle}

\usepackage{amssymb}

\usepackage{amssymb}
\usepackage{amsmath}
\usepackage{color}
\usepackage[nodots]{numcompress}
\newtheorem{theorem}{Theorem}
\newtheorem{definition}{Definition}[section]
\newtheorem{lemma}{Lemma}
\newtheorem{corollary}{Corollary}
\newdefinition{rmk}{Remark}
\newdefinition{example}{Example}
\newproof{pf}{Proof}
\newproof{pot}{Proof of Theorem \ref{thm2}}

\journal{CHAOS, SOLITONS $\&$ FRACTALS}

\begin{document}

\begin{frontmatter}



\title{Co-compact entropy and its variational principle for non-compact spaces}


\author{Zheng Wei\corref{cor1}}
\ead{weizheng@nmsu.edu}
\address{Department of Mathematical Science, New Mexico State University, Las Cruces, New Mexico 88001, U.S.A.}

\author{Yangeng Wang\fnref{fn2,fn3}}
\ead{ygwang@nwn.edu.cn}
\address{Department of Mathematics, Northwest University, Xian, Shaanxi 710069, P.R.China}

\author{Guo Wei\fnref{fn2}}
\ead{guo.wei@uncp.edu}
\address{Department of Mathematics $\&$ Computer Science, University of North Carolina at Pembroke, Pembroke, North Carolina 28372, U.S.A.}

\author{Zhiming Li\fnref{fn2,fn3}}
\ead{lizm@nwu.edu.cn}
\address{Department of Mathematics, Northwest University, Xian, Shaanxi 710069, P.R.China}

\author{Tonghui Wang}
\ead{twang@nmsu.edu}
\address{Department of Mathematical Science, New Mexico State University, Las Cruces,New Mexico 88001, U.S.A.}

\cortext[cor1]{Corresponding author. tel.: (575) 650 6224; fax: (575) 646-1064.}
\fntext[fn2]{The author was supported by the National Natural Science Foundation of China (NSFC No.11371292).}
\fntext[fn3]{The author was supported by the National Natural Science Foundation of China (NSFC No.11301417).}

\begin{abstract}
The purpose of this paper is to generalize the variational principle, which states that the topological entropy is equal to the supremum of the measure theoretical entropies and also the minimum of the metric theoretical entropies, 
to non-compact dynamical systems by utilizing the co-compact entropy. The equivalence between positive co-compact entropy and the existence of $p$-horseshoes over the real line is also proved. This implies a system over real line with positive co-compact entropy contains a compact invariant subset.
\end{abstract}


\end{frontmatter}


\section{Introduction}
In this section, we will go through various definitions and historical notes related to entropy and variational principle. In addition, we review the relationship between the horseshoe and the topological entropy.

\subsection{Entropy and Variatinoal Principle}
In 1965, Adler, Konheim and McAndrew introduced the concept of topological entropy for continuous mappings
defined on compact spaces\cite{Adler}. This entropy is an analogous invariant under conjugation of topological
dynamical systems and can be obtained by maximizing the measure-theoretic entropy over a suitable class of measures
defined on a dynamical system, implying that the topological entropy and measure-theoretic entropy (such as Kolmogorov-Sinai's metric space entropy) are closely related.
Motivated by a conjecture of Adler, Konheim and McAndrew, Goodwyn in 1969 and 1971 compared the
topological entropy  and measure-theoretic entropy  and concluded that topological entropy bounds measure-theoretic entropy\cite{Goodwyn, Goodwyn2}.

In 1970, Dinaburg proved that the topological entropy is equal to the supremum of the measure-theoretic entropy with respect to all invariant measures when the space has finite covering dimension\cite{Dinaburg},  and later in 1970 Goodman  proved the equality in the general case\cite{Goodman}, which is now known as the variational principle on compact space\cite{Walters}. This variational principle plays a fundamental role in ergodic theory and dynamical systems\cite{Ya,Walters}.

In 1971, Bowen considered the entropy for non-compact sets and for group endomorphisms and homogeneous spaces respectively\cite{Bowen5, Bowen2}, which is a generalization of the topological entropy to the non-compact metric spaces, known as a metric space entropy. However, the entropy according to Bowen's definition is metric-dependent and can be positive even for a linear function (Walters' book, pp 169 and 176). In 1973, along with a study of measure-theoretic entropy, Bowen in \cite{Bowen5} gave another definition of topological entropy resembling the Hausdorff dimension, which also equals to the topological entropy defined by Adler, Konheim and McAndrew when the space is compact.

Recently, Canovas and Rodriguez \cite{Canovas}, Malziri and Molaci \cite{Malziri}, Liu, Wang and Wei \cite{Liu}, Wei, Wang and Wei \cite{Wei1,Wei}, and Patr$\tilde{a}$o \cite{Patrao} proposed kinds of definitions of topological entropy on non-compact spaces.

Regarding the variational principle over the non-compact spaces, Patr$\tilde{a}$o extended the variational principle by utilizing the admissible covering and the proper map\cite{Patrao}. In this paper we prove the variational principle for non-compact systems in terms of co-compact entropy,
which is defined by using the perfect map and the co-compact cover. The co-compact entropy was introduced in \cite{Wei1}, and its properties and its relation with other entropies were further investigated in our previous paper \cite{Wei}.

\subsection{Co-compact entropy and horseshoe}
Let $f: I\rightarrow I$ be an interval map, where $I$ is the closed unit interval with the usual Euclidean metric. If $J_1, J_2, \cdots,
J_p$ are disjoint closed non-degenerate subintervals such that
$J_1\cup\cdots\cup J_p\subseteq f(J_i)$ for $i=1,2,\cdots,p$, then
$(J_1, J_2, \cdots, J_p)$ is called a $p$-horseshoe, or simply a
horseshoe if $n= 2$ \cite{S.Ruette}. Block and Coppel called a map with a
horseshoe turbulent \cite{Block}. Misiurewicz showed for a system defined
 on a compact interval, having a positive topological entropy implies that the
system contains the horseshoe phenomenon\cite{M.M2}. On the other hand,
Block, Guckenheimer, Misiurewicz and Young proved that for
a system defined on a compact interval, containing the horseshoe
phenomenon also implies that  the system has a positive topological entropy\cite{Block1}.

Moreover, it is also known that the horseshoe phenomenon implies Devaney chaos\cite{Block}.
However, a system over the whole real line with positive Bowen's entropy
may fail to hold the above horseshoe property. For example, for the
system $(f, d, R)$, where $R$ is the real line with the usual metric
$d(x,y)=|x-y|$ and $f(x)=2x$ for all $x\in R$, its Bowen's entropy
$h_d(f)$ is larger than or equal to $\log 2$ \cite{Walters},  but it
is a simple system without any horseshoe.

It should be pointed out that in the above example, the co-compact entropy (see definition in Section \ref{sec2.2}) is equal to $0$,
as calculated in our previous paper\cite{Wei}. Furthermore, we prove the equivalence between positive co-compact entropy and the existence of $p$-horseshoes over the real line. This also implies a system over real line with positive co-compact entropy contains a compact invariant subset.

\subsection{The purpose, the approach and the outlines}
The main purpose of this paper is to utilize the co-compact covers and obtain a variational principle for non-compact spaces. This result
states that the co-compact entropy is equal to the supremum of the measure theoretic entropies and also the minimum of the metric space entropies. We also show the equivalence between positive co-compact entropy and the existence of $p$-horseshoes over the real line. We proceed by noting that
a perfect map on the real line can be extended to a continuous map on the two-point compactification of $R$ and the set $A=\{-\infty,\infty\}$ is preimage-invariant, i.e., $f^{-1}(A)=A$, and hence the system can be treated in a similar way as in compact space.

Section 2.1 introduces the concept of co-compact open cover of a space and explores the topological properties of such covers.
Section 2.2 introduces the definition of co-compact entropy.
Section 3 shows the variational principle for non-compact spaces under the frame of co-compact entropy.
Sections 4 shows the equivalence between positive co-compact entropy and the existence of $p$-horseshoes over the real. An example is provided to illustrate the result.

\section{Co-compact entropy}         
This section reviews the concept of co-compact entropy and relevant results. The proofs were provided in our previous paper \cite{Wei}.

\subsection{\hspace*{-0.15in}. Co-compact open cover}
Let $(X, f)$ be a topological dynamical system, where $X$ is a Hausdorff space and $f : X \rightarrow X$ is a continuous mapping. We recall the concept of co-compact open covers as follows.

\begin{definition}{\rm\cite{Wei} Let $X$ be a Hausdorff space. For an open subset $U$ of $X$, if $X \backslash U$ is a compact subset of $X$,
then $U$ is called a co-compact open subset. If every element of an open cover $\mathcal{U}$ of $X$  is co-compact,
then $\mathcal {U}$ is called a co-compact open cover of $X$. }
\end{definition}

The following facts are obvious.

(i) The meet of finitely many co-compact open subsets is co-compact, and the union of any collection of co-compact open subsets is
co-compact open\cite{Wei}.

(ii) Let $X$ be Hausdorff. Then any co-compact open cover has a finite subcover\cite{Wei}.

\begin{definition}{\rm\cite{Wei}
Let $X$ and $Y$ be  Hausdorff  spaces and $f : X \rightarrow Y$ a continuous mapping. If $f$ is a closed mapping
and all fibers $f^{-1}(x),  x \in Y$, are compact,  then $f$ is called a perfect mapping. }
\end{definition}

In particular, if $X$ is compact Hausdorff and $Y$ is Hausdorff, every continuous mapping from $X$ into $Y$ is perfect.
If $f:X \rightarrow Y$ is perfect, then $f^{-1}(F)$ is compact for each compact subset $F \subseteq Y$ \cite{Er}.

(iii)
Let $X$ be a Hausdorff space and $f: X \rightarrow X$ a perfect mapping. If $U$ is co-compact  open, then $f^{-1}(U)$ is
co-compact open. Moreover,  If $\mathcal{U}$ is a co-compact open cover of $X$, then  $f^{-1}(\mathcal{U})$ is a co-compact
open cover of $X$\cite{Wei}.

\subsection{The entropy of co-compact open covers}\label{sec2.2}
For compact topological systems $(X, f)$, Adler, Konheim and McAndrew introduced the concept of topological entropy and studied its
properties\cite{Adler}. Their definition is as follows:

For any open cover $\mathcal{U}$ of $X$, denote by
$N_X(\mathcal{U})$ the smallest cardinality of all subcovers of $\mathcal{U}$, i.e.,
$$N_X(\mathcal{U})=\min\{card(\mathcal{V}) : \mathcal{V} {\rm \; is\; a \; subcover\; of\;} \mathcal{U}\}.$$
It is obvious that  $N_X(\mathcal{U})$ is a positive integer. Let $H_X(\mathcal{U})=\log N_X(\mathcal{U})$.
Then $ent(f, \mathcal{U}, X) =$
 $\lim\limits_{n\rightarrow\infty}\frac{1}{n}$ $H_X(\bigvee\limits_{i=0}^{n-1}f^{-i}(\mathcal{U}))$
is called the topological entropy of $f$  relative to
$\mathcal{U}$, and  $ent(f,X)=\sup\limits_\mathcal{U}\{ent(f,\mathcal{U},X)\}$ is called the topological entropy of  $f$.

Now, we generalize Adler, Konheim and McAndrew's entropy to any perfect mappings defined on any Hausdorff space.

Let $X$ be Hausdorff. By the previously stated (ii), when $\mathcal {U}$ is a co-compact open cover of $X$, $\mathcal {U}$
has a finite subcover. Hence, $N_X(\mathcal {U})$, abbreviated as $N(\mathcal {U})$, is a positive integer. Let
$H_X(\mathcal {U}) = \log \: N(\mathcal {U})$, abbreviated as $H(\mathcal {U})$.

For any two given open covers $\mathcal {U}$ and $\mathcal {V}$ of $X$, define

$$\mathcal{U}\bigvee\mathcal{V}=\{U \cap V : U\in\mathcal{U} {\rm \; and\;} V\in\mathcal{V}\}.$$

If for any $U \in \mathcal {U}$, there exists $V \in \mathcal {V}$ such that $U \subseteq V$, then $\mathcal {U}$
is said to be a refinement of $\mathcal {V}${, denoted by} $\mathcal{V}\prec\mathcal{U}$.

The following are some obvious facts:

\smallskip
\smallskip

{\bf Fact 1:} For any open covers $\mathcal{U}$ and $\mathcal{V}$ of $X$,  $\mathcal{U}\prec\mathcal{U}\bigvee\mathcal{V}$.
\smallskip
\smallskip

{\bf Fact 2:} For any open covers $\mathcal{U}$ and $\mathcal {V}$ of $X$,  if $\mathcal{V}$ is a subcover of $\mathcal{U}$, then $\mathcal{U}\prec\mathcal{V}$.
\smallskip
\smallskip

{\bf Fact 3:} For any co-compact open cover $\mathcal{U}$ of $X$,
$H(\mathcal{U})=0 \Longleftrightarrow N(\mathcal{U})=1 \Longleftrightarrow X\in\mathcal{U}$.
\smallskip
\smallskip

{\bf Fact 4:}  For any co-compact open covers $\mathcal{U}$ and $\mathcal {V}$ of $X$,
$\mathcal{V}\prec\mathcal{U}\Rightarrow H(\mathcal{V}) \leq H(\mathcal{U})$.
\smallskip
\smallskip

{\bf  Fact 5:} For any co-compact open covers $\mathcal{U}$ and $\mathcal{V}$,
$H(\mathcal{U}\bigvee\mathcal{V})\leq H(\mathcal{U})+H(\mathcal{V})$.

\smallskip
\smallskip

{\bf Fact 6:} For any co-compact open cover $\mathcal{U}$ of $X$, $H(f^{-1}(\mathcal{U}))\leq H(\mathcal{U})$, and if $f(X) = X$
the equality holds.

\begin{lemma} \label{lemma3.1}
Let $\{a_{n}\}_{n=1}^{ \infty }$ be a sequence of non-negative real numbers satisfying $a_{n+p}\leq a_{n}+a_{p}, n \geq 1, p\geq 1$.
Then the limit $\lim\limits_{n\rightarrow\infty}\frac {a_{n}}{n}$ exists and equals to $\inf \frac {a_{n}}{n}$ (see \rm{\cite{Walters}}).
\end{lemma}

Let $\mathcal {U}$ be a co-compact open cover of $X$. By (iii), for any positive integer $n$ and perfect mapping
$f : X \rightarrow X$, $f^{-1}(\mathcal {U})$ is a co-compact open cover of $X$. On the other hand, $\bigvee\limits_{i=0}^{n-1}f^{-i}(\mathcal{U})$ is a co-compact open cover of $X$. These two facts togehter lead to the following result:

\smallskip
\smallskip

(iv)\label{theorem3.2}
Suppose that  $X$ is Hausdorff. Let $\mathcal{U}$ be a co-compact open cover of $X$, and $f : X \rightarrow X$ a perfect mapping.
Then $\lim\limits_{n\rightarrow\infty}\frac {1}{n}H(\bigvee\limits_{i=0}^{n-1}f^{-i}(\mathcal{U}))$ exists\cite{Wei}.

\bigskip
Next, we introduce the concept of entropy for co-compact open covers.

\begin{definition} {\rm\cite{Wei}
Let $X$ be a Hausdorff space, $f: X\rightarrow X$ be a perfect mapping, and $\mathcal{U}$ be a co-compact open cover of $X$.
The non-negative number $c(f,\mathcal{U})=\lim\limits_{n \rightarrow \infty} \frac {1}{n} H(\bigvee\limits_{i=0}^{n-1}f^{-i}(\mathcal{U}))$ is said to be
the co-compact entropy of $f$ relative to $\mathcal {U}$, and the non-negative number
$c(f)=$ $\sup\limits_{\mathcal {U}}$ $\{c(f, \mathcal{U}) \}$ is said to be the co-compact entropy of $f$, where the supremum is taken over all co-compact open covers.
}
\end{definition}

In particular, when $X$ is compact Hausdorff, any open set of $X$ is co-compact and any continuous mapping $f : X \rightarrow X$
is perfect. Hence, Adler, Konheim and McAndrew's topological entropy is a special case of our co-compact entropy. It should be
aware that the new entropy is well defined for continuous mappings on non-compact spaces, such as  $R^n$, but
Adler, Konheim and McAndrew's topological entropy requires that the space be compact.

Co-compact entropy generalizes the Adler, Konheim and McAndrew's topological entropy,
and yet it holds various similar properties as well, as demostrated by the fact that co-compact entropy
is an invariant of topological conjugation (see (v) below).

Recall that $ent$ denotes Adler, Konheim and McAndrew's topological entropy,
$h_d$ denotes Bowne's entropy and $c$
denotes the co-compact entropy. We summarize some properties of the new entropy from our previous paper\cite{Wei},

\smallskip
\smallskip

(v)\label{theorem3.4}
Let $(X,f)$ and $(Y,g)$  be two topological dynamical systems where
$X$ and $Y$ are Hausodrff, $f : X \rightarrow X$ and $g : Y
\rightarrow  Y$ are perfect mappings. If there exists a
semi-topological conjugation $h: X \rightarrow Y$, then $c(f)\geq
c(g)$. Consequently, when  $h$ is a topological conjugation, we have
$c(f)=c(g)$\cite{Wei}.

\smallskip
\smallskip

(vi)
Let $X$ be Hausdorff and $id : X\rightarrow X$ be the identity mapping. Then $c(id)=0$\cite{Wei}.

When $X$ is Hausdorff and $f:X\rightarrow X$ is perfect, $f^{m}: X \rightarrow X$ is also a perfect mapping \cite{Er}.
\smallskip
\smallskip

(vii)\label{thm1.3.6}
Let $X$ be Hausdorff and $f:X\rightarrow X$ be perfect. Then $c(f^{m})=m\cdot c(f)$\cite{Wei}.
\smallskip
\smallskip

(viii)\label{thm1.3.7}
Let $X$ be Hausdorff and $f:X\rightarrow X$ be perfect. If $\Lambda$ is a closed subset of $X$ and invariant under $f$, i.e.,
$f(\Lambda) \subseteq \Lambda$, then  $c(f |_{\Lambda}) \leq c(f)$\cite{Wei}.

\section{Co-compact entropy and the variational principle}   \label{sec4}
In this section, we will state our main result, the variational principle for the co-compact entropy.

\subsection{Co-compact entropy and its relation with Bowen's entropy}
First let us recall the definition of Bowen's entropy \cite{Bowen6, Walters}. Let $(X, d)$ be a metric space and $f : X \rightarrow X$ a continuous mapping. A compact subset $E$ of $X$ is called a $(n, \epsilon)$-separated set with respect to $f$ if for any different $x, y \in E$, there exists an integer $j$ with $0 \leq j < n$ such that $d(f^j(x), f^j(y)) > \epsilon$. A subset $F$ of $X$ is called a $(n, \epsilon)$-spanning set of a compact set $K$ relative to $f$ if for any $x \in K$, there exists $y \in F$ such that for all $j$ satisfying $0 \leq j < n$,
$d(f^j(x), f^j(y)) \leq \epsilon$.

Let $K$ be a compact subset of $X$. Put
 $$r_{n}(\epsilon,K,f)=\min\{card(F) :  F {\rm \; is\; a\;} (n,\epsilon){\rm -spanning \; set \;for \;} K {\rm \; with\; respect \; to\;} f \},$$
$$s_{n}(\epsilon,K,f)=\max\{card(F) :  F \subseteq K {\rm \; and\;} F {\rm \; is\; a\;} (n,\epsilon) {\rm -separated \; set \; with \; respect \; to \;} f\},$$
$$r(\epsilon,K,f)=\lim\limits_{n\rightarrow\infty}\frac
{1}{n}\log  r_{n}(\epsilon,K,f),~~~~~~~s(\epsilon,K,f)=\lim\limits_{n\rightarrow\infty}\frac {1}{n}\log s_{n}(\epsilon,K,f),$$
$$r(K,f)=\lim\limits_{\epsilon\rightarrow 0} r(\epsilon,K,f) ,~~~~~~~s(K,f)=\lim\limits_{\epsilon\rightarrow 0} s(\epsilon,K,f).$$
Then $\sup\limits_{K}r(K,f)=\sup\limits_{K}s(K,f)$, and this non-negative number denoted by $h_d(f)$ is the Bowen entropy of $f$.

It should be pointed out that Bowen's entropy $h_d(f)$ is metric-dependent, see e.g., \cite{Walters, Liu}. For
the topology of the metrizable space $X$, the selection of different metrics may result in different entropies.

The second definition, which is called $d$-$entropy$ in \cite{Patrao}, is given by
\begin{eqnarray*}
h^d(f)=\sup\limits_{Y}s(Y,f),
\end{eqnarray*}
where now the supremum is taken over all subsets $Y$ of $X$ instead of just the compact ones. Note that the metric $d$ is denoted in the superscript here. Also, $h_d(f)\leq h^d(f)$ always hold and when $X$ is compact, we get the equality. Since $s(Y,f)$ is monotone with respect to $Y$, it follows that
\begin{eqnarray*}
h^d(f)=s(X,f).
\end{eqnarray*}

From our previous paper\cite{Wei}, we have the following result,

(ix) The co-compact entropy is less than or equal to Bowen's entropy.
Precisely, let $(X, d)$ be a metric space, and $f: X\rightarrow X$
be a perfect mapping, then $c(f) \leq  h_{d}(f)$\cite{Wei}.

\subsection{Variational principle for co-compact entropy over non-compact space.}    \label{sec4.2}
From now on we assume that $X$ is a locally compact separable metric space. We denote one-point compactification of $X$ by $\widetilde{X}$, i.e., $\widetilde{X}$ is the disjoint union of $X$ and $\{x_\infty\}$, where $x_\infty$ is some point not in $X$ called the point at infinity. The topology in $\widetilde{X}$ consists of the former open sets in $X$ and the sets $U\cup\{x_\infty\}$, where the complement of $U$ in $X$ is compact.

Let $f:X\rightarrow X$ be a perfect mapping. Define $\widetilde{f}(\widetilde{x}):\widetilde{X}\rightarrow\widetilde{X}$ by
$$
\widetilde{f}(\widetilde{x})=
\begin{cases}
f(\widetilde{x}), & \textrm{if } \widetilde{x}\neq x_\infty, \\
x_\infty, & \textrm{if } \widetilde{x}= x_\infty,
\end{cases}
$$
we have that $\widetilde{f}$ is also a perfect map, called the extension of $f$ to $\widetilde{X}$. In fact we only need to verify that $\widetilde{f}$ is continuous at $x_\infty$. If $U\cup \{x_\infty\}$ is a neighborhood of $x_\infty$, then $X\setminus U$ is compact set in $X$ and we have that
\begin{eqnarray*}
\widetilde{f}^{-1}(U\cup \{x_\infty\})=f^{-1}(U)\cup \{x_\infty\}
\end{eqnarray*}
is also a neighborhood of $x_\infty$, since $f$ is perfect and thus the set $X\setminus f^{-1}(U)$ is also compact. In \cite{Wangsen}, the concept of compact-type metric was introduced. The definition of $\textrm{compact-type metric}$ will be useful:
\begin{definition}\rm\cite{Wangsen}
A metric $d$ of $X$ is of compact-type if it can be extended to a metric $\overline{d}$ of the Alexandroff compactification $\widetilde{X}$. In other words, $d$ is the restriction on $X$ from some metric $\overline{d}$ of $\widetilde{X}$.
\end{definition}

To show the variation principle for co-compact entropy, we need the following lemma:
\begin{lemma} \label{lemma4.3}
Let $d$ be a metric of compact-type on $X$, fix $x_1,x_2,\cdots,x_k\in X$, suppose
$\mathcal{A}_\delta$ $=\{B(x_1,\delta),$ $B(x_2,\delta)$ $,\cdots,$ $B(x_k,\delta)\}$
 is an open cover of $X$ for every $\delta\in (a,b)$, where $0<a<b$, then there exists
  $\delta_{\varepsilon}\in (a,b)$ and an co-compact open cover $\mathcal{A}(\delta_{\varepsilon})$,
   where $\mathcal{A}(\delta_{\varepsilon})$ contains either open ball with radius $\delta_{\varepsilon}$
    or the union of two open balls with radius $\delta_{\varepsilon}$.
  Also, $card(\mathcal{A}_{\delta})=card(\mathcal{A}(\delta_{\varepsilon}))$. Here $B(x,\delta)=\{y\in X|d(x,y)<\delta\}$.
\end{lemma}
\noindent
\pf{
Assume that $\mathcal{A}_\delta=\{B(x_1,\delta),B(x_2,\delta),\cdots,B(x_k,\delta)\}$ is an open cover of $X$, for every $\delta\in (a,b)$, where $0<a<b$. For each fixed $\delta$, the number of balls is finite, it follows that there exists $\delta_{\varepsilon}\in (a,b)$ such that $x_\infty\notin \widetilde{S}(x_{i},\delta_{\varepsilon})$ ( for example, let $\delta_{\varepsilon}$ slightly bigger than $\delta$ ) , for every $i=1,2,\cdots,k,$ where $\widetilde{S}(x_i,\delta_{\varepsilon})=\{\widetilde{y}\in \widetilde{X}|\overline{d}(x_i,\widetilde{y})=\delta_{\varepsilon}\}$. Denoting by $\widetilde{B}(x_i,\delta_{\varepsilon})=\{\widetilde{y}\in \widetilde{X}|\overline{d}(x_i,\widetilde{y})<\delta_{\varepsilon}\}$, there are two cases: (1) the point $x_\infty$ is inside $\widetilde{B}(x_i,\delta_{\varepsilon})$; or (2) the point $x_\infty$ is not in the closure of $\widetilde{B}(x_i,\delta_{\varepsilon})$. In the first case, $X\setminus B(x,\delta_{\varepsilon})=\widetilde{X}\setminus \widetilde{B}(x_i,\delta_{\varepsilon})$, which is compact. Note that there is at least one such open ball, let $B(x^*,\delta_{\varepsilon})$ be any of such ball. In the second case, there exists an open neighborhood $U$ of $x_\infty$ which has empty intersection with $\widetilde{B}(x_i,\delta_{\varepsilon})$. Thus, $B(x,\delta_{\varepsilon})=\widetilde{B}(x_i,\delta_{\varepsilon})$ is in the complement of $U$ in $X$, which has compact closure. Note that the set $B(x,\delta_{\varepsilon})\bigcup B(x^*,\delta_{\varepsilon})$ is co-compact in this case. Now, let
\begin{eqnarray}
\mathcal{A}(\delta_{\varepsilon})&=&\{B(x_i,\delta_{\varepsilon})| x_\infty\in\widetilde{B}(x_i,\delta_{\varepsilon}), B(x_i,\delta_{\varepsilon})\in \mathcal{A}_{\delta_{\varepsilon}}\}\nonumber\\
& &\bigcup\{B(x_i,\delta_{\varepsilon})\cup B(x^*,\delta_{\varepsilon}) |x_\infty\notin\widetilde{B}(x_i,\delta_{\varepsilon}), B(x_i,\delta_{\varepsilon})\in \mathcal{A}_{\delta_{\varepsilon}}\},\nonumber
\end{eqnarray}

then $\mathcal{A}(\delta_{\varepsilon})$ is the co-compact open cover as desired. This completes the proof.
}
\begin{theorem}\label{lemma4.2.1}
Let $f:X\rightarrow X$ be a perfect mapping, let $d$ be a metric of compact-type on $X$. Then it follows that $h^{d}(f)=c(f)$.
\end{theorem}
\noindent
\pf{
We claim that, for all $\varepsilon>0$, there exists $\delta_{\varepsilon}\in (\varepsilon/2,\varepsilon)$ and a co-compact open cover $\mathcal{A}(\delta_{\varepsilon})$ such that $c(f,\mathcal{A}(\delta_{\varepsilon}))\leq r(\varepsilon/2,X,f)$. If $S=\{x_1,x_2,\cdots,x_k\}$ is an $(n,\varepsilon/2)$-spanning set of $X$, for every $\delta\in (\varepsilon/2,\varepsilon)$, we have that
$$\mathcal{B}_\delta=\{B(x_i,\delta)\cap\cdots\cap f^{-n}(B(f^n(x_i),\delta))|x_i\in S\}$$
is a cover of $X$. Also, $\mathcal{A}_\delta=\{B(f^{l}(x_i),\delta))|x_i\in S,0\leq l\leq n\}$ is an open cover of $X$. By Lemma \ref{lemma4.3}, there exists $\delta_{\varepsilon}\in (\varepsilon/2,\varepsilon)$ such that $\mathcal{A}(\delta_{\varepsilon})$ is an co-compact open cover which contains either open ball with radius $\delta_{\varepsilon}$ or the union of two open balls with radius $\delta_{\varepsilon}$. Further more, we have that
$\bigvee\limits_{i=0}^{n}f^{-i}(\mathcal{A}(\delta_{\varepsilon}))\prec \mathcal{B}_{\delta_{\varepsilon}}$, which implies $N(\bigvee\limits_{i=0}^{n}f^{-i}(\mathcal{A}(\delta_{\varepsilon})))\leq card(\mathcal{B}_{\delta_{\varepsilon}})$. Hence for each $(n,\varepsilon/2)-$spanning set $S$ of $X$, there exists $\delta_{\varepsilon}\in (\varepsilon/2,\varepsilon)$ and an co-compact open cover $\mathcal{A}(\delta_{\varepsilon})$ such that $N(\bigvee\limits_{i=0}^{n}f^{-i}(\mathcal{A}(\delta_{\varepsilon})))\leq card(\mathcal{B}_{\delta_{\varepsilon}})\leq card(S)$. Thus it follows that
$N(\bigvee\limits_{i=0}^{n}f^{-i}(\mathcal{A}(\delta_{\varepsilon})))\leq r_n(\varepsilon/2,X,f)$, which implies
\begin{eqnarray}
c(f,\mathcal{A}(\delta_{\varepsilon}))\leq r(\varepsilon/2,X,f).
\end{eqnarray}

Moveover, for the co-compact open cover $\mathcal{A}(\delta_{\varepsilon})$ which has form as defined in Lemma \ref{lemma4.3}, we claim that the inequality $r(|\mathcal{A}_{\delta_{\varepsilon}}|,X,f)$ $\leq c(f,\mathcal{A}(\delta_{\varepsilon}))$ holds, where $|\mathcal{A}_{\delta_{\varepsilon}}|$ is the maximum of the diameters of $B(x_i,\delta_{\varepsilon})\in \mathcal{A}_{\delta_{\varepsilon}}$. In fact, the elements of $\bigvee\limits_{i=0}^{n}f^{-i}(\mathcal{A}(\delta_{\varepsilon}))$ are given by $U_0\cap f^{-1}(U_1)\cap\cdots\cap f^{-n}(U_n)$, where $U_i\in \mathcal{A}(\delta_{\varepsilon})$. From the proof of the Lemma \ref{lemma4.3}, we know each $U_i\in \mathcal{A}(\delta_{\varepsilon})$ has form $U_i=B(x_i,\delta_{\varepsilon})$ or $U_i=B(x_i,\delta_{\varepsilon})\cup B(x^*,\delta_{\varepsilon})$, where $B(x_i,\delta_{\varepsilon}), B(x^*,\delta_{\varepsilon})\in \mathcal{A}_{\delta_{\varepsilon}}.$ Denote $B_{U_i}=B(x_i,\delta_{\varepsilon})$ for each $U_i$.  Now for each $U_0\cap f^{-1}(U_1)\cap\cdots\cap f^{-n}(U_n)\in\bigvee\limits_{i=0}^{n}f^{-i}(\mathcal{A}(\delta_{\varepsilon}))$, let $B_{U_0}\cap f^{-1}(B_{U_1})\cap\cdots\cap f^{-n}(B_{U_n})\in\bigvee\limits_{i=0}^{n}f^{-i}(\mathcal{A}_{\delta_{\varepsilon}})$, such that $B_{U_0}\cap f^{-1}(B_{U_1})\cap\cdots\cap f^{-n}(B_{U_n})\subseteq U_0\cap f^{-1}(U_1)\cap\cdots\cap f^{-n}(U_n)$.
Take an $x\in B_{U_0}\cap f^{-1}(B_{U_1})\cap\cdots\cap f^{-n}(B_{U_n})$ and consider $S$ the set of all such points. We claim that $S$ is an $(n,|\mathcal{A}_{\delta_{\varepsilon}}|)$-spanning set of $X$. In fact, let $y\in X$ and take some $B_{U_0}\cap f^{-1}(B_{U_1})\cap\cdots\cap f^{-n}(B_{U_n})\in\bigvee\limits_{i=0}^{n}f^{-i}(\mathcal{A}_{\delta_{\varepsilon}})$ containing $y$ (since $\mathcal{A}_{\delta_{\varepsilon}}$ itself is a cover). Taking $x\in S$ such that $x\in B_{U_0}\cap f^{-1}(B_{U_1})\cap\cdots\cap f^{-n}(B_{U_n})$, we have that $d(f^i(x),f^i(y))<|\mathcal{A}_{\delta_{\varepsilon}}|$ for every $i=0,1,\cdots, n$. Hence, for the co-compact open cover $\mathcal{A}(\delta_{\varepsilon})$, there exists an $(n,|\mathcal{A}_{\delta_{\varepsilon}}|)$-spanning set $S$ of $X$
such that $r_n(|\mathcal{A}_{\delta_{\varepsilon}}|,X,f)\leq card(S)\leq card(\bigvee\limits_{i=0}^{n}f^{-i}(\mathcal{A}({\delta_{\varepsilon}})))$.
Thus, $r_n(|\mathcal{A}_{\delta_{\varepsilon}}|,X,f)\leq card(\bigvee\limits_{i=0}^{n}f^{-i}(\mathcal{A}({\delta_{\varepsilon}})))$. We get $r_n(|\mathcal{A}_{\delta_{\varepsilon}}|,X,f)\leq N_X(\mathcal{A}(\delta_{\varepsilon}))$.  Taking logarithms, dividing by $n$ and taking limits, it follows as claimed that
\begin{eqnarray}
r(|\mathcal{A}_{\delta_{\varepsilon}}|,X,f)\leq c(f,\mathcal{A}(\delta_{\varepsilon})).
\end{eqnarray}

Since $|\mathcal{A}_{\delta_{\varepsilon}}|\leq 2\delta_{\varepsilon}\leq 2\varepsilon,$ inequality (1) and (2), we have that

$$r(2\varepsilon,X,f)\leq r(|\mathcal{A}_{\delta_{\varepsilon}}|,X,f)\leq c(f,\mathcal{A}(\delta_{\varepsilon}))\leq r(\varepsilon/2,X,f).$$

Taking limits with $\varepsilon\downarrow 0$, it follows that
$$h^d(f)=\lim\limits_{\varepsilon\downarrow 0}c(f,\mathcal{A}(\delta_{\varepsilon}))=\sup\limits_{\varepsilon>0}c(f,\mathcal{A}(\delta_{\varepsilon})).$$

In order to complete the proof, it remains to show that the above supreme is equal to $c(f)$. For any co-compact open cover $\mathcal{A}$ of $X$, by Theorem 4.1 in \cite{Wei}, there exists a Lebesgue number $\varepsilon$ of this cover. We claim that $N(\bigvee\limits_{i=0}^{n}f^{-i}(\mathcal{A}))\leq N(\bigvee\limits_{i=0}^{n}f^{-i}(\mathcal{A}({\delta_{\varepsilon}})))$, where $\mathcal{A}({\delta_{\varepsilon}})$ is a co-compact open cover which has form given in Lemma \ref{lemma4.3}. In fact, since every element of $\bigvee\limits_{i=0}^{n}f^{-i}(\mathcal{A}({\delta_{\varepsilon}}))$ is given by $U_0\cap f^{-1}(U_1)\cap\cdots\cap f^{-n}(U_n)$, we know each $U_i\in \mathcal{A}(\delta_{\varepsilon})$ has form $U_i=B(x_i,\delta_{\varepsilon})$ or $U_i=B(x_i,\delta_{\varepsilon})\cup B(x^*,\delta_{\varepsilon})$. And $B(x_i,\delta_{\varepsilon})$'s are balls of radius $\delta_{\varepsilon}<\varepsilon$, there exists $\{A_0,...,A_n\}\subset\mathcal{A}$ and $A^*\in \mathcal{A}$ such that $B(x_i,\delta_{\varepsilon})\subset A_i$, for each $i$ and $B(x^*,\delta_{\varepsilon})\subset A^*$. Thus we have that
\begin{eqnarray*}
U_0\cap f^{-1}(U_1)\cap\cdots\cap f^{-n}(U_n)&\subseteq & (A_0\cup A^*)\cap f^{-1}(A_1\cup A^*)\cap\cdots\cap f^{-n}(A_n\cup A^*)\\
&\subseteq&(A_0\cap f^{-1}(A_1)\cap\cdots\cap f^{-n}(A_n))\cup (A^*\cap\cdots\cap f^{-n}(A^*))
\end{eqnarray*}
showing that, for each subcover $\mathcal{B}$ of $\bigvee\limits_{i=0}^{n}f^{-i}(\mathcal{A}({\delta_{\varepsilon}}))$, there exists a subcover $\gamma$ of $\bigvee\limits_{i=0}^{n}f^{-i}(\mathcal{A})$ such that $N(\bigvee\limits_{i=0}^{n}f^{-i}(\mathcal{A}))\leq card(\gamma)\leq card(\mathcal{B})$. Hence, we get as claimed that
$N(\bigvee\limits_{i=0}^{n}f^{-i}(\mathcal{A}))\leq N(\bigvee\limits_{i=0}^{n}f^{-i}(\mathcal{A}({\delta_{\varepsilon}})))$. Taking logarithms, dividing by $n$ and taking limits, it follows that
\begin{eqnarray*}
c(f,\mathcal{A})\leq c(f,\mathcal{A}({\delta_{\varepsilon}}))\leq \sup\limits_{\varepsilon>0}c(f,\mathcal{A}(\delta_{\varepsilon})),
\end{eqnarray*}
which shows that
\begin{eqnarray*}
c(f)=\sup\limits_{\varepsilon>0}c(f,\mathcal{A}(\delta_{\varepsilon})).
\end{eqnarray*}
This completes the proof.}

\begin{theorem}\label{theorem4.2.2}
Let $f:X\rightarrow X$ be a perfect mapping, where $X$ is a locally compact separable space. Let $d$ be a metric of compact-type. Then it follows that $ h^{d}(f)=h^{\widetilde{d}}(\widetilde{f})$, where $\widetilde{d}$ is some metric on $\widetilde{X}$, the one-point compactification of $X$, and $\widetilde{f}$ is the extension of $f$ to $\widetilde{X}$. Further more, we have $c(f)=c(\widetilde{f})$.
\end{theorem}
\noindent
\pf{
Since a perfect map is proper, the equality
$h^{d}(f)=h^{\widetilde{d}}(\widetilde{f})$ follows directly from Proposition 2.3 in \cite{Patrao}. And the equality $c(f)=c(\widetilde{f})$ follows from Theorem \ref{lemma4.2.1} and $h^{d}(f)=h^{\widetilde{d}}(\widetilde{f})$.}

\subsection{The variation principle}
In this section we give an extension of the well-known variational principle for non-compact entropy by utilizing co-compact entropy. Unlike topological entropy, metric space entropy doesn't depend on the compactness of the space $X$. As usual, let $(X,\mathcal{B}(X),\mu)$ be a probability space where $\mathcal{B}(X)$ is Borel $\sigma$-algebra and $\mu$ is a $f-$invariant Borel probability measure. For any measurable transformation $f:X\rightarrow X$, we shall denote by $\mathcal{P}(f)$ the collection of all $f$-invariant probability measures on $X$. Given a finite Borel partition $\mathcal{A}$ of $X$,  for every $n\in \mathbb{N}$, define
\begin{eqnarray*}
H(\bigvee_{i=0}^{n-1}f^{-i}(\mathcal{A}))=\sum\limits_{B\in \bigvee_{i=0}^{n-1}f^{-i}(\mathcal{A})}\phi(\mu(B)),
\end{eqnarray*}
where $$\bigvee_{i=0}^{n-1}f^{-i}(\mathcal{A})=\{A_0\cap f^{-1}(A_1)\cap\cdots\cap f^{-n}(A_n)|A_i\in \mathcal{A}\},$$ and $\phi:[0,1]\rightarrow R$ is the continuous function given by,

\begin{equation}
\phi(x)
 = \left\{
\begin{array}{lc}
	-x\log x, & {\rm \; if \; } x\in (0,1]\\
	0,        & {\rm \; if\;} x=0.\\
\end{array}
\right.
\end{equation}
By subadditive ergodic theorem, the following limit
\begin{eqnarray*}
h_\mu(f,\mathcal{A})=\lim\limits_{n\rightarrow \infty}\frac{1}{n}H(\bigvee_{i=0}^{n-1}f^{-i}\mathcal{A})
\end{eqnarray*}
exists.
The metric space entropy of the map $f$ with respect to $\mu$ is thus defined as
\begin{eqnarray*}
h_\mu(f)=\sup\limits_{\mathcal{A}}h_\mu(f,\mathcal{A}),
\end{eqnarray*}
where the supreme is taken over all finite Borel partitions $\mathcal{A}$ of $X$. See \cite{Walters} for more information on metric space entropy. When the space $X$ is not compact, we do not know wether $\mathcal{P}(f)$ is empty or not. However when $f:X\rightarrow X$ is a perfect mapping, and $d$ is a metric of compact-type, by Krylov-Bogolioubov theorem, $\mathcal{P}_{\widetilde{f}}(\widetilde{X})$ is not empty. So by an elegant trick proved by Handel and Kitchens, we have the following result(Lemma 1.5 in \cite{Handel}).

\begin{lemma}\cite{Handel}\label{Lemma4.3.1}
Let $f:X\rightarrow X$ be a perfect mapping such that $\mathcal{P}_{\mu}(f)\neq \emptyset$ and $\widetilde{f}:\widetilde{X}\rightarrow \widetilde{X}$ be its extension in the one-point compactification $\widetilde{X}$ of the locally compact separable space $X$. Then it follows that
\begin{eqnarray*}
\sup\limits_{\mu}h_\mu(f)=\sup\limits_{\widetilde{\mu}}h_{\widetilde{\mu}}(\widetilde{f}),
\end{eqnarray*}
where the supreme are taken, respectively, over $\mathcal{P}_f(X)$ and $\mathcal{P}_{\widetilde{f}}(\widetilde{X})$.
\end{lemma}

Now, we are ready to show the variational principle for locally compact separable space for co-compact entropy. The proof is essentially the same as  Theorem 3.2 in \cite{Patrao}, the main difference lies on the ways to define \lq\lq non-compact\rq\rq  topological entropy and the different assumptions on the map $f$.
\begin{theorem}\label{thmvp}
Let $f:X\rightarrow X$ be a perfect mapping, where $X$ is a locally compact separable space. Then it follows that
\begin{eqnarray*}
\sup\limits_{\mu}h_\mu(f)=c(f)=\min\limits_{d}h_{d}(f),
\end{eqnarray*}
where the minimum is attained whenever $d$ is a compact-type metric.
\end{theorem}
\noindent
\pf{ The variational principle for compact spaces states that
\begin{eqnarray*}
c(\widetilde{f})=\sup\limits_{\widetilde{\mu}}h_{\widetilde{\mu}}(\widetilde{f}).
\end{eqnarray*}
By Proposition 1.4 in \cite{Handel}, we have that
\begin{eqnarray*}
\sup\limits_{\mu}h_\mu(f)\leq\inf\limits_{d}h_{d}(f).
\end{eqnarray*}
Applying Lemma \ref{Lemma4.3.1}, it follows that
\begin{eqnarray*}
c(\widetilde{f})=\sup\limits_{\mu}h_\mu(f)\leq\inf\limits_{d}h_{d}(f).
\end{eqnarray*}
Applying Theorem \ref{lemma4.2.1} and Theorem \ref{theorem4.2.2}, we have
\begin{eqnarray*}
h_d(f)\leq h^d(f)=c(f)=\sup\limits_{\mu}h_\mu(f)\leq\inf\limits_{d}h_{d}(f),
\end{eqnarray*}
where, in the  second term, $d$ is any compact-type metric.}

\begin{rmk}  On locally compact separable space, Patr$\tilde{a}$o defined the admissible entropy by using admissible cover, which is a finite cover contains not only co-compact open sets but also open sets with compact closure(note that the co-compact cover does not require finiteness and only contains co-compact open sets)\cite{Patrao}. Since the admissible entropy also admits the variational principle(Theorem 3.2 in \cite{Patrao}), we conclude that the co-compact entropy equals to the admissible entropy in this case. As it will be shown in the next section, a system with positive co-compact entropy is equivalent to the existence of the $p$-horseshoe. However, to the best of our knowledge, we haven't seen any result about the relation between admissible entropy and $p$-horseshoe in the literature.
\end{rmk}

\section{Positive co-compact entropy implies horseshoes over the real line}
For maps defined on closed intervals, the results by Misiurewicz\cite{M.M2} and
Block-Guckenheimer- Misiurewicz -Young\cite{Block1} disclose a relation between
positive entropy and the existence of p-horseshoes. However, these
results are unavailable for interval maps defined on the whole real
line (Adler, Konheim and McAndrew's topological entropy is only
defined for compact systems and, for Bowen's metric space entropy, there
are counterexamples with positive entropies but without $p$- horseshoes). 
To generalize these results to interval maps defined on the real line,
the co-compact entropy will be utilized in this section. It will be proved that for a perfect mapping $f$ from the
real line into itself, if it has a positive co-compact entropy, then
it has horseshoes; conversely, if it has a p-horseshoe, then its
co-compact entropy is larger than or equal to $\log p$.

\subsection{Entropy and horseshoe}
The following two Theorems show the equivalence between the positive entropy and the existence of p-horseshoes over compact interval,

\begin{theorem}\label{Misiurewiczthm} (Misiurewicz's Theorem \cite{M.M2})
Let $f: I\rightarrow I$ be an interval map of positive entropy
$h(f)$ (Adler's topological entropy or, which is the same, Bowen's
metric space entropy). Then for every $\lambda $ satisfying $0<\lambda
<h(f)$ and every  positive integer $N$, there exists disjoint closed
intervals $J_{1}, ..., J_{p}$ and a  positive integer $n\geq N$,
such that $(J_{1}, ..., J_{p})$ is a $p-$horseshoe for $f^{n}$ and
$\frac{\log p}{n}\geq \lambda$.
\end{theorem}

\begin{theorem}\label{Misiurewiczthm1}  
(The special case of Block-Guckenheimer-Misiurewicz-Young's
Theorem \cite{Block1}) Let $f: I\rightarrow I$ be an
interval map. If $f$ has a $p$-horseshoe,  then the entropy (Adler's
topological entropy or, which is the same, Bowen's metric space entropy)
$h(f)$ is larger than or equal to $\log p$.
\end{theorem}
Moreover, it is also known that the horseshoe phenomenon implies
chaos \cite{Block1}. However, a system over the whole real line with positive Bowen's entropy may fail to hold the above horseshoe property. Motivated by this, we prove the equivalence between the positive co-compact entropy
and the existence of p-horseshoes over the real line in the following section.

\subsection{Co-compact entropy and horseshoe}
The purpose of this section is to establish two results: Generalize the results stated in the
above two theorems from maps defined on closed intervals
to maps defined on the real line by
utilizing the concept co-compact entropy.


Let $f:R\rightarrow R$ be a perfect mapping. By the property of perfect mappings, a map $f: R \to R$ is perfect implies $f$ converges at infinity.
There are 4 different cases to consider:
\smallskip
\smallskip

(C1). $\lim\limits_{x\rightarrow \pm\infty}f(x)=-\infty$,
\smallskip
\smallskip

(C2). $\lim\limits_{x\rightarrow \pm\infty}f(x)=+\infty$,
\smallskip
\smallskip

(C3). $\lim\limits_{x\rightarrow -\infty}f(x)=-\infty$ and $\lim\limits_{x\rightarrow +\infty}f(x)=+\infty$,
\smallskip
\smallskip

(C4). $\lim\limits_{x\rightarrow -\infty}f(x)=+\infty$ and $\lim\limits_{x\rightarrow +\infty}f(x)=-\infty$.
\smallskip
\smallskip

Define $h:R\rightarrow (0,1)$ by $h(x)=\frac{1}{\pi}\tan^{-1}(x)+\frac{1}{2}$. $h$ is a homeomorphism
between $R$ and $(0,1)$, set $g=h\circ f\circ h^{-1}$, then the system $(R,f)$ is topological conjugate to the system $((0,1),g)$. Note that $f$ converges at infinity implies $g$
 converges at $0$ and $1$. Therefore, we have corresponding 4 different situations for $g$:
\smallskip
\smallskip

(A1). $\lim\limits_{x\rightarrow 0}g(x)=0$ and $\lim\limits_{x\rightarrow 1}g(x)=0$,
\smallskip
\smallskip

(A2). $\lim\limits_{x\rightarrow 0}g(x)=1$ and $\lim\limits_{x\rightarrow 1}g(x)=1$,
\smallskip
\smallskip

(A3). $\lim\limits_{x\rightarrow 0}g(x)=0$ and $\lim\limits_{x\rightarrow 1}g(x)=1$,
\smallskip
\smallskip

(A4). $\lim\limits_{x\rightarrow 0}g(x)=1$ and $\lim\limits_{x\rightarrow 1}g(x)=0$.
\smallskip
\smallskip

Let's define the extension $\widetilde{g}(x):[0,1]\rightarrow [0,1]$ of $g$ by
$$
\widetilde{g}(x)=
\begin{cases}
g(x), & \textrm{if } x\neq 0 \textrm{ and } 1,\\
\lim\limits_{x\rightarrow 1}g(x), & \textrm{if } x=1,\\
\lim\limits_{x\rightarrow 0}g(x), & \textrm{if } x=0.\\
\end{cases}
$$
Also, $g$ converges at $0$ and $1$ implies that $\widetilde{g}$ is continuous at $0$ and $1$ and the set $A=\{0,1\}$ is
preimage-invariant under $\widetilde{g}$, i.e., $\widetilde{g}^{-1}(A)=A$.

Theorem \ref{thm5.1} below generalizes Theorem \ref{Misiurewiczthm}.

\begin{theorem}\label{thm5.1}
Let $f:R\rightarrow R$ be a perfect mapping with co-compact entropy $c(f)>0$. Then for every $\lambda $ satisfying $0<\lambda
<c(f)$ and every  positive integer $N$, there exists disjoint closed
intervals $J_{1}, ..., J_{p}$ and a  positive integer $n\geq N$,
such that $(J_{1}, ..., J_{p})$ is a $p-$horseshoe for $f^{n}$ and
$\frac{log p}{n}\geq \lambda$.
\end{theorem}
\noindent
\pf{ Since the co-compact entropy and $p$-horseshoe are invariant under topological conjugation(see \cite{Wei,S.Ruette}), and
 $(R,f)$ is topological conjugate to the system $((0,1),g)$, where $g$ is defined above, it is enough to show the result for
the system $((0,1),g)$.

I claim that $c(\widetilde{g})\geq c(g)$. In fact, let $\mathcal{U}$ be any co-compact open cover of $(0,1)$.
For each $U\in\mathcal{U}$, $U\cup \{0, 1\}$ is an open set in $[0,1]$.
Define $\mathcal{U}^{*}=\{U\cup \{0,1\}| U\in \mathcal{U}\}$. $\mathcal{U}^{*}$ is an open cover for
$[0,1]$ and $Card(\mathcal{U})= Card(\mathcal{U}^*)$. For any $n$, we also have
$N(\bigvee\limits_{i=0}^{n-1}g^{-i}(\mathcal{U}))= N(\bigvee\limits_{i=0}^{n-1}\widetilde{g}^{-i}(\mathcal{U}^*))$(since $A$
is $\widetilde{g}$ invariant). Taking logarithms, dividing by $n$ and taking limits, it follows that
$c(g,\mathcal{U})=c(\widetilde{g},\mathcal{U}^*)$.
 Therefore, for each co-compact open cover $\mathcal{U}$ of $(0,1)$, there exists an open cover $\mathcal{U}^*$
  of $[0,1]$ such that $c(g,\mathcal{U})=c(\widetilde{g},\mathcal{U}^*)$, it follows as claimed that
   $c(\widetilde{g})=\sup\{c(\widetilde{g},\mathcal{U}^*)\}\geq \sup\{c(g,\mathcal{U})\}= c(g)$.

Now, we know $c(\widetilde{g})\geq c(g)>0$. By applying Theorem \ref{Misiurewiczthm} for $\widetilde{g}$ on
$[0,1]$, we get for every $\lambda $ satisfying $0<\lambda
<c(\widetilde{g})$ and every  positive integer $N$, there exists disjoint closed
intervals $J_{1}, ..., J_{p}$ in $[0,1]$ and a  positive integer $n\geq N$,
such that $(J_{1}, ..., J_{p})$ is a $p-$horseshoe for $\widetilde{g}^{n}$ and
$\frac{log p}{n}\geq \lambda$. Also, since the set $A$ is preimage-invariant under $\widetilde{g}$, this implies none of $p$ disjoint closed
intervals $J_{1}, ..., J_{p}$ contains $0$ or $1$. Therefore, $J_{i}\subset (0,1)$ and this completes the proof.}\\

\noindent
\noindent

Let $\sum=\{0,1\}^{\mathbb{N}}$ define the shift map $\sigma:\Sigma\rightarrow\Sigma$ by $\sigma((a_{n})_{n\in \mathbb{N}})=(a_{n+1})_{n\in \mathbb{N}}$. The dynamical
system $(\Sigma,\sigma)$ is a shift on two letters. Block showed that if an interval map has a horseshoe formed of two disjoint subintervals, then there exists an invariant compact subset on which the action of the map is almost a shift on two letters\cite{Block,Block2}. Therefore, the following result is a straightforward consequence of Theorem \ref{thm5.1},
\begin{corollary}
Let $f:R\rightarrow R$ be a perfect mapping with co-compact entropy $c(f)>0$.
Then there exists a compact invariant subset of $R$ with respect to $f$.
\end{corollary}

Theorem \ref{thm5.2} below generalizes Theorem \ref{Misiurewiczthm1}. We need the following definition and Lemma.
Let $f:R\rightarrow R$ be a perfect map, if $J_{0}, J_{1}, ...,
J_{n-1}$ are subintervals on $R$ such that $J_{i}\subset
f(J_{i-1})$, for $1\leq i \leq n-1$, then $(J_{0},
 J_{1}, ..., J_{n-1})$ is called a chain of intervals on $R$\cite{S.Ruette}.
Notice that in the proof of Lemma \ref{Lemma3.3}, the assumption that
$f$ is a perfect map has been used.

\begin{lemma} \label{Lemma3.3}
Let $f:R\rightarrow R$ be a perfect map, $(J_{0}, J_{1}, ...,
J_{n-1})$ is a chain of intervals on $R$. Then there is a
subinterval $K\subset J_{0}$, such that $f^{n-1}(K)=J_{n-1}$ and
$f^{i}(K)\subset J_{i}$ for all $0\leq i \leq n-1$.
\end{lemma}

\noindent
\pf{
First we prove that this lemma holds for
$n=2$.

Case (1): $J_{0}$ and $J_{1}$ are both bounded. See proof of Lemma
1 in \cite{M.M2}.

Case (2): $J_{0}$ is a bounded interval and $J_{1}$ is an unbounded
interval. Note here the function $f$ is defined over $R$,
$f(\overline{J_{0}})$ will be a bounded interval. This case is
impossible.

Case (3): $J_{0}$ is an unbounded interval and $J_{1}$ is a bounded
interval. In this case, $J_{0}$ could be $(-\infty,v)$,
$(-\infty,v]$, $[u,\infty)$, $(u,+\infty)$, $(-\infty,+\infty)$. and
$J_{1}$ could be $(a,b)$, $[a,b)$, $(a,b]$, $[a,b]$.

Assume $J_{0}=(-\infty,v)$. When $J_{1}=[a,b]$, since $J_{1}\subset
f(J_{0})$, there exists $x,y\in (-\infty,v)$ such that $f(x)=a,
f(y)=b$. Then the proof is same as case 1. When $J_{1}=(a,b]$, if
$a$ belongs to $f(J_{0})$, then we go back to the previous proof.
Now assume $a$ doesn't belongs to $f(J_{0})$. Since $J_{1}\subset
f(J_{0})$, there exists $y\in (-\infty,v)$ such that $f(y)=b$. Then
define $y'=\min\{(-\infty,y]\cap f^{-1}(b)\}$ (note that we can take
min because $f$ is a perfect map and $f^{-1}(b)$ is a compact set),
$y''=\sup\{[y,v)\cap f^{-1}(b)\}$ (note if $y''=v$, we only need to
define $y'$). Set $K=(-\infty,y'], K'=[y'',v)$, then one of $K$ and
$K'$ will do the job (since only boundary points of $J_{0}$ can
approach to $a$). Same proof works for the case when $J_{1}=[a,b)$
and $J_{1}=(a,b)$.

The cases when $J_{0}=(-\infty,v]$, $(-\infty,v]$,
$(-\infty,+\infty)$, $(u,+\infty)$ and $[u,+\infty)$ follow a
similar proof as $J_{0}=(-\infty,v)$.

Case(4): both $J_{0}$ and $J_{1}$ are unbounded intervals. In this
case, $J_{0}$ could be $(-\infty,v)$, $(-\infty,v], $ $[u,\infty),$
$(u,+\infty), $ $(-\infty,+\infty)$ and $J_{1}$ could be
$(-\infty,b)$, $(-\infty,b]$, $[a,\infty)$, $(a,+\infty)$,
$(-\infty,+\infty)$.

Assume $J_{0}=(-\infty,v)$. When $J_{1}=(-\infty,b]$, since
$J_{1}\subset f(J_{0})$, there exists $y\in (-\infty,v)$ such that
$f(y)=b$. Then define $y'=\min\{(-\infty,y]\cap f^{-1}(b)\}$ (note
that we can take min because $f^{-1}(b)$ is a compact set). Set
$K=(-\infty,y']$, then $K$ will do the job. When
$J_{1}=(-\infty,b)$, if $b$ belongs to $f(J_{0})$, it goes back to
previous case. Now assume $b$ doesn't belong to $f(J_{0})$, then
$f((-\infty,v))=J_{1}$. Same proof works when $J_{1}=[a,\infty),
(a,+\infty)$. When $J_{1}=(-\infty,+\infty)$, then
$f((-\infty,v))=J_{1}$.

The cases $J_{0}=(-\infty,v]$, $[u,\infty), $ $(u,+\infty)$,
$(-\infty,+\infty)$ follow a similar proof of the case
$J_{0}=(-\infty,v)$.

Suppose that the Lemma is true for $n$ and consider a chain of
intervals $(J_{0},..., J_{n-1}, J_{n})$. By assumption, there exists
an interval $K\subset J_{0}$ with $f^{n-1}(K) = J_{n-1}$ and
$f^{i}(K)\subset J_{i}$ for $0\leq i \leq n-1$, thus one has
$f^{n}(K)\supset J_{n}$. Applying the proof $n=2$ to $f^{n}$ and the
chain of interval $(K,J_{n})$, we obtain an interval $K'\subset K$
such that $f^{n}(K')=J_{n}$. This completes the proof.

\begin{theorem} \label{thm5.2}
Let $f:R\rightarrow R$ be
a perfect map. If $f^n$ has a $p-$horseshoe for some $n\geq 1$, then
$c(f)\geq \frac{\log p}{n}$.
\end{theorem}

\noindent
\pf{
Put $g=f^{n}$. Then $g$ is still a perfect
map. Let $J_{1},\cdots,J_{p}$ be a $p-$horseshoe for $g$, i.e.
$J_{1}, J_{2},
 \cdots, J_{p}$ are $p$ disjoint closed and bounded intervals
 such that $J_{1} \cup J_{2} \cup \cdots \cup J_{p}
\subset g(J_{i})$ for $1\leq i \leq p$. Let $U_{1}=(R \setminus
J_{2}) \cap (R \setminus J_{3})\cdots \cap (R \setminus J_{p}),$
$U_{2}=(R \setminus J_{1}) \cap (R \setminus J_{3}) \cdots \cap (R
\setminus J_{p}),$ $\cdots$, $U_{p}=(R \setminus J_{1}) \cap (R
\setminus J_{2}) \cdots \cap (R \setminus J_{p-1}).$ Since $J_{1},
J_{2}, \cdots, J_{p}$ are disjoint closed intervals,
$\mathcal{U}=\{U_{1}, U_{2}, \cdots, U_{p}\}$ is a co-compact open
cover for $R$. Also for each $J_{i}$, there is a unique element
$U_{i}\in\mathcal{U}$ such that $J_{i}\subset U_{i}$. For every
$(i_{0},\cdots,i_{n-1})\in \{ 1,\cdots,p\}^{n}$, define
$$J_{i_{0},\cdots,i_{n-1}} = \{\ x\in R \mid g^{k}(x)\in J_{i_{k}},
0\leq k \leq n-1 \}.$$

Since $J_{1},\cdots,J_{p}$ is a $p-$horseshoe, thus
$J_{i_{0}},\cdots,J_{i_{n-1}}$ is a chain of interval. By Lemma
\ref{Lemma3.3}, there exists $J$ such that $f^{j}(J)\subset
J_{i_{j}}$ for $0\leq j\leq n-1$ and $f^{n-1}(J)= J_{i_{n-1}}$.
Therefore, $J\subset J_{i_{0},\cdots, i_{n-1}}$ is not empty.
Moreover $J_{i_{0},\cdots, i_{n-1}}$ is contained in a unique
element of $\bigvee\limits_{i=0}^{n-1}g^{-i}(\mathcal{U})$, which is
$U_{i_{0}} \cap g^{-1}(U_{i_{1}}) \cap\cdots\cap
g^{-(n-1)}(U_{i_{n-1}})$. Hence $N_{n}(\mathcal{U}) \geq p^{n}$ for
all integers $n$, and $c(g)\geq c(g,\mathcal{U})=\lim\limits_{n\rightarrow
+\infty} \frac{\log N_{n}(\mathcal{U})}{n}\geq \log p$. Recall property (vi) of the co-compact entropy, $c(f)=\frac{c(f^n)}{n}\geq \frac{\log p}{n}$.
This completes the proof.}

\begin{example}
Let $f: R \longrightarrow R$ be a perfect map defined by each $n\in N$:

$$
f(x) =
\begin{cases}
3x-2n,           &\text{ if } x\in [n, n+\frac{1}{3}]\\
-3x+4n+2,               &\text{ if } x\in (n+\frac{1}{3}, n+\frac{2}{3}]\\
3x-2n-2,                &\text{ if }  x\in (n+\frac{2}{3}, n+1].
\end{cases}
$$
\end{example}

Clearly, $[n, n+1]$ is a closed invariant set of $f$ for every $n$.
Consider the restriction of $f$ on the unit interval $[n, n+1]$. It
is known that the topological entropy of $f|_{[n, n+1]}$ is at least
$\log 2$ (See \cite{S.Ruette}). Since the space $[n, n+1]$ is compact, the
co-compact entropy is equal to the topological entropy, we have
$c(f|_{[0, 1]}) \geq \log 2$.  By Property (viii), $c(f)\geq
c(f |_{[0,1]}) \geq \log2$. By Theorem \ref{thm5.1}, $f^{m}$ contains
horseshoes for some positive integer $m$. In fact,
$[n,n+\frac{1}{3}], [n+\frac{2}{3}, n+1]$ is a $2$-horseshoes for
$f$.

\section{Concluding remarks}
In the application of dynamical systems, spaces that are often encountered are
non-compact manifolds such as $R$. Consequently, finding an appropriate entropy and
developing a variational principle for these manifolds become important and useful.
For non-compact systems, the co-compact entropy retains various fundamental properties
of the topological entropy; in particular, it only depends on the topology, not on
the selection of a metric when the space is metrizable. A noticeable property of the co-compact entropy
is that, as proved in this paper, it admits the variational principle for non-compact manifolds, which states that
the co-compact entropy is equal to the supremum of the measure theoretical entropies
and also the minimum of the metric space entropies, where the latter entropy is metric-dependent.

In addition, the co-compact entropy provides a numerical characteristic/sign for
the existence of horseshoes. This implies that there exists
a compact invariant subset for a system with positive co-compact entropy. These advantages
demonstrate the appropriateness of using the co-compact entropy in studying chaotic phenomena over
non-compact spaces. As a future study, the relation between the co-compact entropy and
other chaotic phenomena (e.g., Li-York chaos, Devaney chaos, and bifurcations
of chaotic scattering etc.) over non-compact spaces need to be explored.






\begin{thebibliography}{99}
\bibliographystyle{model3-num-names}
\bibliography{<your-bib-database>}




\bibitem{Adler} Adler RL, Konheim AG, McAndrew MH. Topological entropy. Trans. Amer. Math. Soc. 1965;114:309-319.

\bibitem{Goodwyn} Goodwyn LW. Topological entropy bounds measure-theoretic entropy. Proc Amer. Math. Soc. 1969;23:679-688.

\bibitem{Goodwyn2} Goodwyn LW. Comparing topological entropy with measure-theoretic entropy. Amer. J Math 1972;94:366-388.

\bibitem{Dinaburg} Dinaburg EI. The relation between topological entropy and metric entropy. Soviet Math. 1970;11:13-16.

\bibitem{Goodman} Goodman TNT. Relating topological entropy and measure entropy. Bull. Lond. Math. Soc. 1971;3:176-180.

\bibitem{Walters} Walters P. An introduction to ergodic theory. Berlin:Springer; 1982.

\bibitem{Ya} Pesin YaB. Dimension Theory in Dynamical Systems, Contemporary Views and Applications. Chicago:University of Chicago Press; 1997.

\bibitem{Bowen5} Bowen R. Topolog entropy for non-compact sets. Trans. Amer. Math. Soc. 1971;14:401-414.

\bibitem{Bowen2} Bowen R. Entropy follr group endomorphisms and homogeneous spaces. Trans. Amer. Math. Soc. 1971;14:401-14.

\bibitem{Canovas} Canovas JS, Rodriguez JM. Topological entropy of maps on the real line. Top Appl. 2005;153:735-746.

\bibitem{Malziri} Malziri M, Molaci MR. An extension of the notion of the topological entropy. Chaos, Solitions $\&$ Fractals 2008;36:370–373.

\bibitem{Liu} Liu L, Wang Y,  Wei G. Topological entropy of continuous functions on topological spaces. Chaos, Solitions $\&$ Fractals 2009;39:417–427.

\bibitem{Wei1} Wei Z, Wang Y, Wei G. A new entropy: co-compact entropy. Science Journal of Northwest University Online 2009;7:1-13.

\bibitem{Patrao} Patrao M. Entropy and its variational principle for non-compact metric spaces. Ergod. Th. $\&$ Dynam. Sys. 2010;30:1529-1542.

\bibitem{S.Ruette} Ruette S. Chaos for continuous interval maps, http://www.math.u-psud.fr/~ruette/, 2002.

\bibitem{Block} Block LS, Coppel WA.  Dynamics in one dimension, Lecture Notes in Mathematics vol 1513. Berlin: Springer; 1992.

\bibitem{M.M2} Misiurewicz M. Horseshoes for continuous mappings of an interval. Bull. Acad. Polish Sci. 1979;27:167-169.

\bibitem{Block1} Block LS, Guckenheimer J, Misiurewicz M, Young L.S. Global Theory of Dynamical Systems, Lecture Notes in Mathematics. No. 819. New York: Springer-Verlag; 1980.

\bibitem{Wei} Wei Z, Wang Y,  Wei G, Wang T,  Bourquin S. The entropy of co-compact opencovers, Entropy 2013;15:2464-2479.

\bibitem{Er} Engelking R. General Topology. Warszawa:PWN; 1989.

\bibitem{Bowen6} Bowen R. Topological entropy and Axiom A , Global Analysis. Proc. Sympos. Pure Math., Amer. Math. Soc. 1970;14:23-42.

\bibitem{Wangsen}  Wang Y, Wei G, Campbell WH. Sensitive dependence on initial
 conditions between dynamical systems and their induced hyperspace dynamical systems. Top. Appl. 2009;156:803-811.

\bibitem{Handel} Handel M, Kitchens B. Metrics and entropy for non-compact spaces. Israel J. Math. 1995;91:253-271.

\bibitem{Block2} Block LS. Homoclinic points of mappings of the interval. Proc. Amer. Math. Soc. 1978;72:576-580.
\end{thebibliography}
\end{document}